%% file: frac_qxfu.tex
\documentclass[12pt]{article}
\usepackage{amsmath,amssymb,amsthm,a4wide,color,graphicx}
\newtheorem{theorem}{Theorem}[section]
\newtheorem{lemma}{Lemma}[section]
\newtheorem{corollary}{Corollary}[section]

\newtheorem{remark}{Remark}[section]
\def\du{\underline{du}}
\def\dq{\underline{dq}}
\def\df{\underline{df}}
\def\ulu{\underline{u}}
\def\olu{\overline{u}}

\title{Reconstruction of space-dependence and nonlinearity of a reaction term in a subdiffusion equation}
\author{Barbara Kaltenbacher\footnote{
Department of Mathematics,
Alpen-Adria-Universit\"at Klagenfurt.
barbara.kaltenbacher@aau.at.}
\and
William Rundell\footnote{
Department of Mathematics,
Texas A\& M University, College Station, Texas 77843 USA.
rundell@math.tamu.edu}
}
\begin{document}
\maketitle
\begin{abstract}
In this paper we study the simultaneous reconstruction of two coefficients in a reaction-subdiffusion equation, namely a nonlinearity and a space dependent factor. The fact that these are coupled in a multiplicative matter makes the reconstruction particularly challenging.
Several situations of overposed data are considered: boundary observations over a time interval, interior observations at final time, as well as a combination thereof.
We devise fixed point schemes and also describe application of a frozen Newton method. In the final time data case we prove convergence of the fixed point scheme as well as uniqueness of both coefficients.
Numerical experiments illustrate performance of the reconstruction methods, in particular  dependence on the differentiation order in the subdiffusion equation.
\end{abstract}

\input{intro}

\subsection{The inverse problem}\label{sec:intro}
Our aim is to identify $q(x)$ and $f(u)$ in  
\begin{equation}\label{PDE}
\begin{aligned}
\partial_t^\alpha u-\triangle u +q(x) f(u) &= r \quad&& \text{ in }\Omega\times(0,T)\\
Bu&=0&&\text{ on }\partial\Omega\times(0,T)\\
u(t=0)&=u_0&& \text{ in }\Omega
\end{aligned}
\end{equation}
from (a) time trace (b) final time observations (c) mixed observations 
\begin{equation}\label{obs}
\begin{aligned}
&(a) \ h^i(x,t)=u^i(x,t), \ x\in \Gamma\subseteq\partial\Omega, \ t\in(0,T),\\
&(b) \ g^i(x) = u^i(x,T), \ x\in \Omega, \\
&(c) \ h(x,t) = u^1(x,t), \ x\in \Gamma\subseteq\partial\Omega, \ t\in(0,T),\qquad
g(x) = u^1(x,T), \ x\in \Omega
\end{aligned}
\end{equation}
with two excitation pairs $(r^1,u_0^1)$, $(r^2,u_0^2)$ (in case of (c) we only need one of them).
In \eqref{PDE}, the boundary operator $B$ defines Dirichlet, Neumann or impedance boundary conditions
\[Bv=\partial_\nu v+\gamma v,\] 
with 
$\gamma$ a non-negative measurable mapping defined on $\partial\Omega$
and the Dirichlet case 
corresponding to $\gamma=\infty$. 

Identifying both $q$ and $f$ can only be unique up to a constant factor: If the pair $(q(x),f(u))$ solves this inverse problem then so does the pair $(cq(x),\frac{1}{c}f(u))$ for any $c\in\mathbb{R}\setminus\{0\}$.
This ambiguity can be resolved by some sort of normalisation; we do so by fixing the value of $q$ at some point $\bar{x}$ as 
\begin{equation}\label{qbar}
q(\bar{x})=\bar{q} 
\end{equation}
with some prescribed value $\bar{q}$, for example obtained from some direct observation at a boundary point $\bar{x}$.

The forward operator $F=(F^1,F^2)$ is thus defined by 
\begin{equation}\label{F}
F(q,f)= C ((u^i)_{i\in I}) =\begin{cases}
(\text{tr}_{\Gamma\times(0,T)} u^1, \ \text{tr}_{\Gamma\times(0,T)} u^2)&(a)\\
(\text{tr}_{\Omega\times\{T\}} u^1, \ \text{tr}_{\Omega\times\{T\}} u^2)&(b)\\
(\text{tr}_{\Gamma\times(0,T)} u^1, \ \text{tr}_{\Omega\times\{T\}} u^1)&(c)
\end{cases}
\end{equation}
where $u^i$ solves \eqref{PDE} with $r=r^i$, $u_0=u_0^i$, $I=\{1,2\}$ in case (a) or (b), $I=\{1\}$ in case (c).
We denote the data by $y$ in each of these cases
\begin{equation}\label{y}
y= \begin{cases} (h^1,h^2)&(a)\\ (g^1,g^2)&(b)\\ (h^1,g^1)&(c)
\end{cases}
\end{equation}

\section{Analysis of the forward problem}\label{sec:analysis}
In this section we will prove well-definedness and Lipschitz continuous Fr\'{e}chet differentiability of the forward operator. The latter can be viewed as a basis for applying Newton's method to the inverse problem, which we will consider in Section~\ref{sec:Newton}. 

To establish well-posedness of the nonlinear initial boundary value problem \eqref{PDE}, we will make use of \cite[Theorems 2, 4]{KubickaYamamoto2018} for the linear problem
\begin{equation}\label{PDE_lin}
\begin{aligned}
\partial_t^\alpha (u-u_0)-\triangle u +c(x,t)u &= r \quad&& \text{ in }\Omega\times(0,T)\\
Bu&=0&&\text{ on }\partial\Omega\times(0,T)\\
u(t=0)&=u_0&& \text{ in }\Omega
\end{aligned}
\end{equation}
which we here summarize for the convenience of the reader. It is formulated for homogeneous Dirichlet boundary conditions $B=\text{tr}_{\partial\Omega}$, $H^1_B(\Omega)=H_0^1(\Omega)$ there, but extends to homogeneous Neumann or impedance boundary conditions in a straightforward way.
To take these conditions into account, we denote by $-\triangle_B$ the weak form of the negative Laplacian under these boundary conditions  
\[
\langle -\triangle_B v,w\rangle_{(H^1_B)^*,H^1_B} = \int_\Omega \nabla v\cdot\nabla w\, dx 
+\int_{
\{0<\gamma<\infty\}}\gamma vw\ dS \text{ for all }w\in H^1_B,
\]
where the superscript ${}^*$ denotes the topological dual and
\[
H^1_B=\{w\in H^1(\Omega)\, : \, w=0\text{ on }
\{\gamma=\infty\}\},
\]
and we assume $\gamma\vert_{\{0<\gamma<\infty\}}\in L^\infty(\{0<\gamma<\infty\})$.
We here use the notation $u\in {}_{u_0}H^\alpha(0,T;X)$ for $I_t^{1-\alpha}(u-u_0)\in {}_0H^1(0,T;Z)=\{v\in H^1(0,T;Z)\, : \, v(t=0)=0\}$.
The low and high regularity regularity results from \cite{KubickaYamamoto2018} are submerged into the following theorem as cases (lo) and (hi).

\begin{theorem}\label{thm:KY17}(\cite[Theorems 2, 4]{KubickaYamamoto2018}; particularised to $L=-\triangle + c\,\text{id}$)\hfill\break\\
Assume that $\alpha\in(0,1)$, $T>0$, $\Omega\subseteq\mathbb{R}^d$, $\Omega$ Lipschitz, 
$p\in(\frac{d}{2},\infty]$ and $c\in L^\infty(0,T;L^p(\Omega))$.
\begin{enumerate}
\item[(lo)] If $u_0\in L^2(\Omega)$, $r\in W_{lo}:=L^2(0,T;H^1_B(\Omega)^*)$, 
then there exists a unique weak solution 
\[
u\in U_{lo}:={}_{u_0}H^\alpha(0,T;H^1_B(\Omega)^*)\cap H^{\alpha/2}(0,T;L^2(\Omega))\cap L^2(0,T;H^1_B(\Omega))
\]
of \eqref{PDE_lin} and $\|u\|_{U_{lo}}\leq C(\|u_0\|_{L^2(\Omega)}+\|r\|_{L^2(0,T;H^1_B(\Omega)^*)})$
holds for some $C$ independent of $u_0$ and $r$.
\item[(hi)] If additionally $\Omega\in C^{1,1}$, 
$p\geq2$, 
$u_0\in H^1_B(\Omega)$, $r\in W_{hi}:=L^2(0,T;L^2(\Omega))$,
then \[
u\in U_{hi}:={}_{u_0}H^\alpha(0,T;L^2(\Omega))\cap H^{\alpha/2}(0,T;H^1_B(\Omega))\cap L^2(0,T;H^2(\Omega))
\]
and $\|u\|_{U_{hi}}\leq C(\|u_0\|_{H^1_B(\Omega)}+\|r\|_{L^2(0,T;L^2(\Omega))})$
holds for some $C$ independent of $u_0$ and $r$.
\end{enumerate}
\end{theorem}

\medskip

In order to show well-posedness of the nonlinear problem \eqref{PDE}, we apply the Implicit Function Theorem 
to the mapping
\[
H:X\times U\to W, \quad (q,f,u) \mapsto \partial_t^\alpha (u-u_0)-\triangle_B u + q(x) f(u) - r
\]
which implicitly defines the parameter-to-state map $S:X\mapsto U$ by 
\[
H(q,f,S(q,f)) = 0.
\]
The function spaces are chosen as 
\begin{equation}\label{X}
X=L^p(\Omega)\times W^{1,\infty}_{\textrm{aff}}(\mathbb{R}), \quad (U,W)\in\{(U_{lo},W_{lo}),\,(U_{hi},W_{hi})\},
\end{equation}
with $p$, $U_{lo}$, $U_{hi}$, $W_{lo}$, $W_{hi}$ as in Theorem~\ref{thm:KY17}. 
Here $U^0$ is defined like $U$, but with $u_0$ replaced by $0$ and $W^{1,\infty}_{\textrm{aff}}(\mathbb{R})$ includes affine functions $f(\zeta)=f_0+\mathfrak{s}\,\zeta$ for some $f_0\in\mathbb{R}$, $\mathfrak{s}\in\mathbb{R}$
\begin{equation}\label{UW2inftyaff}
\begin{aligned}
&U^0_{lo}:={}_{0}H^\alpha(0,T;H^1_B(\Omega)^*)\cap H^{\alpha/2}(0,T;L^2(\Omega))\cap L^2(0,T;H^1_B(\Omega))\\
&U^0_{hi}:={}_{0}H^\alpha(0,T;L^2(\Omega))\cap H^{\alpha/2}(0,T;H^1_B(\Omega))\cap L^2(0,T;H^2(\Omega))\\
&W^{k,\infty}_{\textrm{aff}}(\mathbb{R}):=\{f\text{ measurable }\, ; \, \|f\|_{W^{k,\infty}_{\textrm{aff}}(\mathbb{R})}:=
|f(0)|+
\|f'\|_{W^{k-1,\infty}(\mathbb{R})}
<\infty\}, \quad k\in\mathbb{N}
\end{aligned}
\end{equation}
Note that a function $f\in W^{k,\infty}_{\textrm{aff}}$ is not necessarily essentially bounded, but satisfies the growth estimate
\begin{equation}\label{estfval}
|f(\zeta)|\leq|f(0)|+\|f'\|_{L^\infty(\mathbb{R})}|\zeta|\leq \|f\|_{W^{1,\infty}_{\textrm{aff}}(\mathbb{R})}(1+|\zeta|)
\end{equation}

\medskip

Well-definedness of $H:X\times U\to W$ follows by bounding the nonlinear term  
\begin{equation}\label{qfu}
\|qf(u)\|_{W}\leq
\|qf(0)\|_W+\|f'\|_{L^\infty(\mathbb{R})}\|q\,u\|_W
\end{equation}
and further estimating as follows.\hfill\break\\
In the low regularity scenario $U=U_{lo}$, $W=W_{lo}$, with $p_2\in[1,\infty]$ such that the embedding $H^1(\Omega)\to L^{p_2}(\Omega)$ (and therewith by duality, with $p_2^*=\frac{p_2}{p_2-1}$, the embedding $L^{p_2^*}(\Omega)\to H^1(\Omega)^*$, $p_2^*\geq\frac{2d}{d+2}$, $p_2\leq\frac{2d}{d-2}$, with strict inequality in case $d=2$) is continuous we use 
\begin{equation}\label{estp3}
\begin{aligned}
&\|q\, u\|_{L^2(0,T;H^1(\Omega)^*)}\lesssim \|q\, u\|_{L^2(0,T;L^{p_2^*}(\Omega))}
\leq \|q\|_{L^p(\Omega)} \|u\|_{L^2(0,T;L^{p_3/2}(\Omega))}\\
&\|q\, u\, v\|_{L^2(0,T;H^1(\Omega)^*)}\lesssim \|q\, u\, v\|_{L^2(0,T;L^{p_2^*}(\Omega))}
\leq \|q\|_{L^p(\Omega)} \|u\|_{L^4(0,T;L^{p_3}(\Omega))}
\|v\|_{L^4(0,T;L^{p_3}(\Omega))}, \\ 
&p_3=\frac{2pp_2^*}{p-p_2^*},
\end{aligned}
\end{equation}
where by interpolation $H^{\alpha/2}(0,T;L^2(\Omega))\cap L^2(0,T;H^1_B(\Omega))\subseteq
H^{\theta\alpha/2}(0,T;H^{1-\theta}(\Omega))\subseteq L^4(0,T;L^{p_3}(\Omega))$ with $\theta=\frac{1}{2\alpha}$, $1-\frac{1}{2\alpha}-\frac{d}{2}\geq-\frac{d}{p_3}$, that is 
\begin{equation}\label{p_lo}
\frac{1}{p}\leq\frac{3}{d}-\frac{1}{d\alpha}-\frac12.
\end{equation}
In the limiting case $\alpha\to1$ this allows for $p\geq6$ if $d=3$ and $p>2$ if $d=2$.
\hfill\break\\ 
In the high regularity scenario $U=U_{hi}$, $W=W_{hi}$, we use
\begin{equation}\label{estp4}
\begin{aligned}
&\|q\, u\|_{L^2(0,T;L^2(\Omega))}
\leq \|q\|_{L^p(\Omega)} \|u\|_{L^2(0,T;L^{p_4/2}(\Omega))}\\
&\|q\, u\, v\|_{L^2(0,T;L^2(\Omega))}
\leq \|q\|_{L^p(\Omega)} \|u\|_{L^4(0,T;L^{p_4}(\Omega))}
\|v\|_{L^4(0,T;L^{p_4}(\Omega))}, \\ 
&p_4=\frac{4p}{p-2},
\end{aligned}
\end{equation}
where by interpolation $H^{\alpha/2}(0,T;H^1_B(\Omega))\cap L^2(0,T;H^2(\Omega))\supseteq
H^{\theta\alpha/2}(0,T;H^{2-\theta}(\Omega))\supseteq L^4(0,T;L^{p_4}(\Omega))$ with $\theta=\frac{1}{2\alpha}$, $2-\frac{1}{2\alpha}-\frac{d}{2}\geq-\frac{d}{p_4}$, that is 
\begin{equation}\label{p_hi}
\frac{1}{p}\leq\frac{4}{d}-\frac{1}{d\alpha}-\frac12.
\end{equation}
Thus, as compared to \eqref{p_lo} we can use a smaller integrability index $p$ as compared to the low regularity scenario; this comes at the cost of needing more regularity on $u_0$ and $r$, cf. Theorem~\ref{thm:KY17}. 
In particular, in the limiting case $\alpha\to1$, $p\geq2$ if $d=3$ and actually $p=2$ can be used in case $d=2$ provided $\alpha$ is close enough to one (recall that $p\geq2$ was imposed in part (hi) of Theorem~\ref{thm:KY17} anyway). 

\medskip

To verify the assumptions of the Implicit Function Theorem, we first of all establish existence of a reference point $(q^0,f^0,u^0)\in X\times U$ such that $H(q^0,f^0,u^0)=0$. 
We can do so in several alternative ways: 
\begin{enumerate}
\item[(i)] by choosing $q^0\in L^p(\Omega)$, $\mathfrak{s}\in\mathbb{R}$ arbitrarily, setting $f^0(\zeta)=\mathfrak{s}\cdot \zeta$ and using Theorem~\ref{thm:KY17} with $c(x,t)=\mathfrak{s}\,q^0(x)$ to obtain existence of $u^0\in U$ such that $H(q^0,f^0,u^0)=0$;
\item[(ii')] by assuming that a solution $(q_{\textrm{act}},f_{\textrm{act}},u_{\textrm{act}})$ to the inverse problem under consideration exists and setting $(q^0,f^0,u^0)=(q_{\textrm{act}},f_{\textrm{act}},u_{\textrm{act}})$;
\item[(ii)] by slightly relaxing (ii'), assuming existence of $(q^0,f^0,u^\sharp)$ such that $\|H(q^0,f^0,u^\sharp)\|_{W}$ is small enough
and applying the Newton-Kantorowich Theorem to the mapping  $G_{(q^0,f^0)}:U\to W$, $u\mapsto G_{(q^0,f^0)}(u):= H(q^0,f^0,u)$, while using the fact (see below) that $G_{(q^0,f^0)}'(u^\sharp)=\frac{\partial H}{\partial u}(q^0,f^0,u^\sharp):U^0\to W$ is an isomorphism.
\end{enumerate}
Secondly, the linear operator $(q^0,f^0,u^\sharp)\in X\times U$  
\[
\frac{\partial H}{\partial u}(q^0,f^0,u^\sharp):U^0\to W, \quad 
\du \mapsto \partial_t^\alpha \du -\triangle_B \du + q^0 {f^0}'(u^\sharp) \du
\] 
is an isomorphism due to Theorem~\ref{thm:KY17} with $c(x,t)=q^0(x){f^0}'(u^\sharp(x,t))$ which due to $q^0\in L^p(\Omega)$, ${f^0}'(u^\sharp)\in L^\infty(0,T;L^\infty(\Omega))$ indeed defines an element $c\in L^\infty(0,T;L^p(\Omega))$.
\hfill\break\\
Thirdly, to establish continuity (even Lipschitz continuity) of the derivatives, 
we estimate
\[
\begin{aligned}
&\|\frac{\partial H}{\partial u}(q^1,f^1,u^1)\du-\frac{\partial H}{\partial u}(q^2,f^2,u^2)\du\|_W\\
&\leq\|(q^1-q^2) {f^1}'(u^1)\du\|_W + \|q^2({f^1}-{f^2})'(u^1)\du\|_W + \|q^2 ({f^2}'(u^1)-{f^2}'(u^2)) \du\|_W\\
&\leq L \|(q^1,f^1,u^1)-(q^2,f^2,u^2)\|_{X\times U}\|\du\|_U,
\end{aligned}
\]
since we have 
\begin{equation}\label{estfLinfty}
\begin{aligned}
&\|{f^1}'(u^1)\|_{L^\infty(0,T;L^\infty(\Omega))} \leq \|{f^1}'\|_{L^\infty(\mathbb{R})}
\\
&\|({f^1}-{f^2})'(u^1)\|_{L^\infty(0,T;L^\infty(\Omega))} \leq \|(f^1-f^2)'\|_{L^\infty(\mathbb{R})}\\
&\|{f^2}'(u^1)-{f^2}'(u^2)\|_{L^P(0,T;L^Q(\Omega))} \leq \|{f^2}''\|_{L^{\infty}(\mathbb{R})} \|u^1-u^2\|_{L^P(0,T;L^Q(\Omega))},
\end{aligned}
\end{equation}
together with \eqref{estp3}, \eqref{estp4}.
\hfill\break\\
Boundedness and Lipschitz continuity of the other partial derivatives
\[
\begin{aligned}
&\frac{\partial H}{\partial q}(q^0,f^0,u^\sharp):L^p(\Omega)\to W, \quad 
\dq \mapsto \dq f^0(u^0), \\
&\frac{\partial H}{\partial f}(q^0,f^0,u^\sharp):W^{1,\infty}(\mathbb{R})\to W, \quad 
\df \mapsto q^0 \df(u^0)
\end{aligned}
\] 
follows from the estimates
\[
\begin{aligned}
&\|\frac{\partial H}{\partial q}(q^1,f^1,u^1)\dq-\frac{\partial H}{\partial q}(q^2,f^2,u^2)\dq\|_W
=\|\dq (f^1(u^1)-f^2(u^2)\|_W\\
&\leq\|\dq({f^1}-{f^2})(u^1)\|_W + \|\dq ({f^2}(u^1)-{f^2}(u^2))\|_W
\end{aligned}
\]
\[
\begin{aligned}
&\|\frac{\partial H}{\partial f}(q^1,f^1,u^1)\df-\frac{\partial H}{\partial f}(q^2,f^2,u^2)\df\|_W
=\|q^1 \df(u^1) - q^2 \df(u^2)\|_W\\
&\leq\|(q^1-q^2)\df(u^1)\|_W +  \|q^2 (\df(u^1)-\df(u^2))\|_W
\end{aligned}
\]
combined with \eqref{qfu}, 
as well as \eqref{estp3} or \eqref{estp4}.
 
Thus we can apply the Implicit Function Theorem, which also yields Lipschitz continuity of $S'$, due to 
\[
\begin{aligned}
S'(q^1,f^1)-S'(q^1,f^1)
=&-{K_1}^{-1}\frac{\partial H}{\partial (q,f)}(q^1,f^1,u^1)
+{K_2}^{-1}\frac{\partial H}{\partial (q,f)}(q^2,f^2,u^2)\\
=&-{K_1}^{-1}\left(\frac{\partial H}{\partial (q,f)}(q^1,f^1,u^1)
-\frac{\partial H}{\partial (q,f)}(q^2,f^2,u^2)\right)\\
&+{K_2}^{-1}(K_1-K_2){K_1}^{-1}\frac{\partial H}{\partial (q,f)}(q^2,f^2,u^2)\\
\text{ with }&K_i=\frac{\partial H}{\partial u}(q^i,f^i,u^i).
\end{aligned}
\]

In fact, as revealed by \eqref{estfLinfty}, higher order regularity ${f^2}''\in L^\infty(\mathbb{R})$ is only required at the reference point that we can fix to $f^2=f^0$, as 
the proof of the Implicit Function Theorem shows: For any 
$(q,f)\in \mathcal{B}_{\rho_X}^X(q^0,f^0)\subseteq X$, $u^1,u^2\in \mathcal{B}_{\rho_U}^U(u^0)\subseteq U$, 
the fixed point operator defined by  
\[
T_{q,f}:U\to U, \quad T_{q,f}(u)=u-\frac{\partial H}{\partial u}(q^0,f^0,u^0)^{-1}H(q,f,u)
\]
satisfies
\[\begin{aligned}
&T_{q,f}(u^1)-T_{q,f}(u^2)\\
&=-\frac{\partial H}{\partial u}(q^0,f^0,u^0)^{-1}\bigl(H(q,f,u^1)-H(q,f,u^2)-\frac{\partial H}{\partial u}(q^0,f^0,u^0)(u^1-u^2)\bigr)\\
&=-\frac{\partial H}{\partial u}(q^0,f^0,u^0)^{-1}\bigl((H(q,f,u^1)-H(q,f,u^2))-(H(q^0,f^0,u^1)-H(q^0,f^0,u^2))\\
&\hspace*{4cm}+\,
H(q^0,f^0,u^1)-H(q^0,f^0,u^2)-\frac{\partial H}{\partial u}(q^0,f^0,u^0)(u^1-u^2)\bigr),
\end{aligned}\]
thus, by the Mean Value Theorem 
applied to the mappings 
\[\begin{aligned}
&\tau\mapsto H(q^0+\tau(q-q^0),f^0+\tau(f-f^0),u^1)-H(q^0+\tau(q-q^0),f^0+\tau(f-f^0),u^2),\\ 
&\tau\mapsto H(q^0,f^0,u^1+\tau(u^2-u^1))-\frac{\partial H}{\partial u}(q^0,f^0,u^0)(u^1+\tau(u^2-u^1))
\end{aligned}\]
the estimate
\[\begin{aligned}
&\|T_{q,f}(u^1)-T_{q,f}(u^2)\|_U \leq\|\frac{\partial H}{\partial u}(q^0,f^0,u^0)^{-1}\|_{L(U,W)}\\
&\times\Bigl(
\sup_{\tau\in[0,1]}\|\frac{\partial H}{\partial(q,f)} (q^0+\tau(q-q^0),f^0+\tau(f-f^0),u^1)\\
&\hspace*{2cm}
-\frac{\partial H}{\partial(q,f)} (q^0+\tau(q-q^0),f^0+\tau(f-f^0),u^2)\|_{L(X,W)}\|(q-q^0,f-f^0)\|_X\\
&\hspace*{1cm}
+\sup_{\tau\in[0,1]}\|\frac{\partial H}{\partial u} (q^0,f^0,u^1+\tau(u^2-u^1))-\frac{\partial H}{\partial u} (q^0,f^0,u^0)\|_{L(U,W)}\|u^1-u^2\|_U \Bigr)\\
&\leq C(\|(q-q^0,f-f^0)\|_X+\|u^1-u^0\|_U+\|u^2-u^0\|_U)\, \|u^1-u^2\|_U
\end{aligned}\]
according to the Lipschitz estimates above, with a constant depending only on $\|(q^0,f^0)\|_X$, $\|u^0\|_U$ and $\|{f^0}''\|_{L^\infty(\mathbb{R})}$. With $\rho_X\geq\|(q-q^0,f-f^0)\|_X$, $\rho_U\geq\|u^1-u^0\|_U,\|u^2-u^0\|_U$ small enough, this yields contractivity of $T_{q,f}$ as well as its invariance on $\mathcal{B}_{\rho_X}^X(q^0,f^0)$ and therefore existence and uniqueness of a fixed point $u=S(q,f)\in\mathcal{B}_{\rho_U}^U(u^0)$.


\begin{theorem}\label{thm:S}
Assume that $q^0\in L^p(\Omega)$, $f^0\in W^2_{\textrm{aff}}(\mathbb{R})$ and either
\begin{itemize}
\item[(i)] $f^0$ affine, that is, $f^0(\zeta)=f_0+\mathfrak{s}\,\zeta$ for some $f_0\in\mathbb{R}$, $\mathfrak{s}\in\mathbb{R}$ or
\item[(ii)] $\|\frac{\partial H}{\partial u}(q^0,f^0,u^\sharp)^{-1}\|\, L\, \|\frac{\partial H}{\partial u}(q^0,f^0,u^\sharp)^{-1}H(q^0,f^0,u^\sharp)\|_{X\times U}\leq\frac12$ for some $u^\sharp\in U$.
\end{itemize}
Then there exists $\rho_X>0$ such that the parameter-to-state map 
\[
S:\mathcal{B}_{\rho_X}^X(q^0,f^0)(\subseteq X)\to U, \quad (q,f)\mapsto u \text{ solving \eqref{PDE} }
\]
is well defined and Lipschitz continuously Fr\'{e}chet differentiable in either of the following two scenarios
\begin{itemize}
\item[(lo)] $U=U_{lo}$, $r\in W_{lo}$, $p$ in the definition of $X$ satisfies \eqref{p_lo};
\item[(hi)] $U=U_{hi}$, $r\in W_{hi}$, $p$ in the definition of $X$ satisfies \eqref{p_hi}.
\end{itemize}
\end{theorem}
\begin{remark}\label{rem:f}
Note that the condition $f'\in W^{1,\infty}(\mathbb{R})$ imposes a bound on the growth of $f'$; however, Theorem~\ref{thm:S} can still be applied in the context of the common polynomial choices of $f$ by using a smooth cutoff function $\eta_M\in C^2(\mathbb{R})$ such that $\eta(u)=1$ for $|u|\leq M$, $\eta(u)=0$ for $|u|\geq M+1$ for some $M$ large enough, $\|\eta\|_{L^\infty(\mathbb{R})}\leq1$ and setting   
$f_M(\zeta)= f(0)+f'(0)\zeta+\int_0^\zeta \int_0^\xi \eta_M(v) f''(v)\, dv$, so that 
$f_M'(\zeta)= f'(\zeta)$ for $|\zeta|\leq M$ and $f_M'(\zeta)=C_M:=f'(0)+\int_0^{M+1} \eta_M(v) f''(v)\, dv$ for $|u|\geq M+1$, $\|f_M''\|_{L^\infty(\mathbb{R})}\leq\|f''\|_{L^\infty(0,M+1)}$.

Alternatively to this growth bound, monotonicity conditions on $f$ could be used to establish well-posedness of the forward problem. As in the integer derivative case $\alpha=1$, \eqref{PDE} is well-posed if $q\geq0$ and $f$ is monotonically inceasing. Without going into detail about the proof, we just give an intuition by deriving the relevant energy estimate that results from testing with $u$ (assuming $f(0)=0$ for simplicity of exposition):
\[
\begin{aligned}
\int_0^T\int_\Omega r\, u\, dx\, dt
=\int_0^T\Bigl(\int_\Omega \Bigl(\partial_t^\alpha u\, u +|\nabla u|^2+q \,f(u)\,u \, dx\Bigr)\, dx 
+\int_{
\{0<\gamma<\infty\}}\gamma \, u^2\ dS\Bigr)\, dt\\
\geq \cos(\pi\alpha/2)\|u\|_{H^{\alpha/2}(0,T;L^2(\Omega))}^2+\|\nabla u\|_{L^2(0,T;L^2(\Omega))}^2
\end{aligned}
\]
due to coercivity of the Abel integral operator \cite[Lemma 2.3]{Eggermont:1987}; see also \cite[Theorem 1]{VoegeliNedaiaslSauter:2016}
and monotonicity $f(u)\,u=(f(u)-f(0))\,(u-0)\geq0$.
\end{remark}

\section{Iterative reconstruction schemes}\label{sec:iterations}
While Newton's method, which we briefly discuss in Section~\ref{sec:Newton}, is a generic approach for nonlinear operator equations, we will also derive some particularly tailored fixed point schemes in Section~\ref{sec:fixedpoint}. The principle behind these schemes is to project the PDE onto the observation manifold and extract a fixed point equation for the searched for quantities $q$ and $f$ from this. In case (b) of final time data, the equation resulting from projection is still nonlinear and we apply the logarithm to write it as an equivalent linear relation. This is not needed in case (c) of mixed data, which allows for partial elimination. 
Note that in case (a) of time trace observations, this projection approach is not applicable.    
  
\subsection{Frozen Newton method}\label{sec:Newton}
A Newton step for solving the operator equation formulation $F(q,f)=y$ of the inverse problem with $F$, $y$ as in \eqref{F}, \eqref{y}  is defined by updating 
\begin{equation}\label{Newton}
\begin{aligned}
&q^{(k+1)}=q^{(k)}+s_q^{(k)}, \quad f^{(k+1)}=f^{(k)}+s_f^{(k)}, \quad \\
&\text{ where } F'(q^{(k)},f^{(k)})(s_q^{(k)},s_f^{(k)})=y-F(q^{(k)},f^{(k)}).
\end{aligned}
\end{equation}
Linearising $F$ at some $(q^0,f^0)$ amounts to evaluating $F'(q^0,f^0)(\dq,\df)=C ((u^i)_{i\in I})$, where $\du^i$ solves 
\begin{equation}\label{dPDE}
\begin{aligned}
\partial_t^\alpha \du^i-\triangle \du + q^0(x){f^0}'(u^{0i})\du^i&= 
-\dq(x) f^0(u^{0i}) - q^0(x) \df(u^{0i})\quad&& \text{ in }\Omega\times(0,T)\\
B\du&=0&&\text{ on }\partial\Omega\times(0,T)\\
\du(t=0)&=0&& \text{ in }\Omega
\end{aligned}
\end{equation}
and $u^{0i}$ solves \eqref{PDE} with $r=r^i$, $u_0=u_0^i$, $q=q^0$, $f=f^0$.
With the special choice $f^0(u)=u$, both \eqref{dPDE} and \eqref{PDE} become linear
\begin{equation}\label{dPDE_qcfu}
\begin{aligned}
\partial_t^\alpha u^{0i} +  (-\triangle +q^0)u^{0i} &= r^i\\
\partial_t^\alpha \du^i + (-\triangle +q^0)\du^i&= -\dq(x)\, u^{0i} - q^0\, \df(u^{0i}).
\end{aligned}
\end{equation}
Setting up the linear operator $F'(q^{(k)},f^{(k)})$ thus involves the (numerical) solution of \eqref{dPDE} with $(q^0,f^0)=(q^{(k)},f^{(k)})$ for each direction $(\dq,\df)\in X$. To avoid this substantial computational effort, instead of \eqref{Newton}
we fix the evaluation point of $F'$ and implement a frozen Newton method
\begin{equation}\label{frozenNewton}
\begin{aligned}
&q^{(k+1)}=q^{(k)}+s_q^{(k)}, \quad f^{(k+1)}=f^{(k)}+s_f^{(k)}, \quad \\
&\text{ where } F'(q^0,f^0)(s_q^{(k)},s_f^{(k)})=y-F(q^{(k)},f^{(k)}).
\end{aligned}
\end{equation}
Since the linear operator equation to be solved in \eqref{frozenNewton} is still ill-posed, each step of \eqref{frozenNewton} needs to be equipped with some regularization and the overall iteration with a criterion for early stopping to avoid unbounded propagation of noise, see e.g., \cite{BakKok04,KNS08} and the references therein.
A convergence proof of such iterative methods in the sense of regularization cf. \cite{EnglHankeNeubauer:1996,KNS08} requires restrictions on the nonlinearity of the operator such as the so-called tangential cone condition or range invariance of the linearisation.
However, the problem under consideration seems to be too nonlinear to satisfy any of these, unless full observations $C=\text{id}_U$ are taken.

\subsection{Fixed point schemes}\label{sec:fixedpoint}
\subsubsection*{case (a) time trace data}
Projecting the PDE onto the observation manifold $\Gamma\times(0,T)$ in the time trace case 
\begin{equation}\label{projPDE_timetrace}
\begin{aligned}
&-\partial_t^\alpha h^{1}(t)+\triangle u^1(x_0,t) +r^1(x_0,t) = q(x_0) f(h^1(t))\quad && x_0\in\Gamma, \ t\in(0,T)\\
&-\partial_t^\alpha h^{2}(t)+\triangle u^2(x_0,t) +r^2(x_0,t) = q(x_0) f(h^2(t))\quad && x_0\in\Gamma, \ t\in(0,T)
\end{aligned}
\end{equation}
does not lead to a system that allows to extract $q(x)$ for all $x\in\Omega$.
So we are restricted to using Newton's method in that situation.

\subsubsection*{case (b) final time data}
Fixed point schemes in this setting have been developed in \cite{ZhangZhi:2017,frac_potential} for the identification of $q$ only and  \cite{frac_reactiondiffusion,frac_reacdiff_sys} for recovering $f$ only.

Projecting the PDE onto the observation manifold yields
\begin{equation}\label{projPDE_finaltime}
\begin{aligned}
&-\partial_t^\alpha u^1(x,T)+\triangle g^1(x)+r^1(x,T) = q(x) f(g^1(x))&& x\in\Omega\\
&-\partial_t^\alpha u^2(x,T)+\triangle g^2(x)+r^2(x,T) = q(x) f(g^2(x))&& x\in\Omega.
\end{aligned}
\end{equation}
To extract $q$, $f$ from the left hand side, we take the logarithm of the absolute values on both sides
\[
\begin{aligned}
&\tilde{q}(x)+\tilde{f}(g^1(x))=\tilde{r}^{1,(k)}(x)\\
&\tilde{q}(x)+\tilde{f}(g^2(x))=\tilde{r}^{2,(k)}(x)\\
&\text{where }
\tilde{q}=\log\circ q, \quad \tilde{f}=\log\circ f, \\ 
&\hspace*{1.2cm}\tilde{r}^i(x)=\log(|-\partial_t^\alpha u^i(x,T)+\triangle g^i(x)+r^i(x,T)|), \quad i\in\{1,2\}
\end{aligned}
\]
and subtract, which eliminates $q(x)$.
We thus consider the fixed point scheme
\begin{equation}\label{fixedpoint_finaltime}
\begin{aligned}
&\tilde{f}^{(k+1)}(g^1(x))-\tilde{f}^{(k+1)}(g^2(x))=\tilde{r}^{1,(k)}(x)-\tilde{r}^{2,(k)}(x)\\
&\tilde{q}^{(k+1)}(x)= \tilde{r}^{1,(k)}(x) - \tilde{f}^{(k+1)}(g^1(x))\text{ or }
\tilde{q}^{(k+1)}(x)= \check{\tilde{r}}^{1,(k)}(x) - \tilde{f}^{(k+1)}(g^1(x))\\
&\hspace*{1cm}\text{where }
\tilde{r}^i(x)=\log(|-\partial_t^\alpha u^{i,(k)}(x,T)+\triangle g^i(x)+r^i(x,T)|), \quad i\in\{1,2\},\\
&\hspace*{2.2cm}
\check{\tilde{r}}^1(x)=\log(|-\partial_t^\alpha \check{u}^{1,(k)}(x,T)+\triangle g^1(x)+r^1(x,T)|).
\end{aligned}
\end{equation}
Here $u^{i,(k)}$ solves \eqref{PDE} with $q=\exp(\tilde{q}^{(k)})$, $f=\exp(\tilde{f}^{(k)})$ and $\check{u}^{1,(k)}$ solves \eqref{PDE} with $q=\exp(\tilde{q}^{(k)})$, $f=\exp(\tilde{f}^{(k+1)})$; that is, we may choose to update the PDE solution before updating $q$. 

Alternatively, we consider the overall system  
\begin{equation}\label{fixedpoint_finaltime_overall}
\begin{aligned}
&\tilde{q}^{(k+1)}(x)+\tilde{f}^{(k+1)}(g^i(x))=\tilde{r}^{i,(k)}(x), \quad i\in\{1,2\}
\end{aligned}
\end{equation}
We can write this in an incremental way, using the (projected) PDE, 
to arrive at
\begin{equation}\label{fixedpoint_finaltime_overall_incr}
\begin{aligned}
&\tilde{q}^{(k+1)}(x)+\tilde{f}^{(k+1)}(g^i(x))\\
&=\log(|\exp(\tilde{q}^{(k)}(x)) \exp(\tilde{f}^{(k)}(u^i(x,T))) + \triangle (g^i(x,T)-u^i(x,T))|)
, \quad i\in\{1,2\}
\end{aligned}
\end{equation}
which reveals the fact that the iteration is driven by the data misfit $g^i(x,T)-u^i(x,T)$

To prove that the first equation in \eqref{fixedpoint_finaltime} determines $\tilde{f}^{(k+1)}$, we assume the value $f_0$ of $f$ to be known at some point that is contained in the range of $g^1$
\[
f(s_0):=f_0>0 \ \text{ with } [\ulu,\olu]:=g^1(\Omega), \quad s_0\in [\ulu,\olu], 
\] 
restrict ourselves to the spatially 1-d case $\Omega=(0,1)$ and assume that $g^1$, $g^2$ $\in C^1(\Omega)$ are strictly monotone and that there exist constants $0<c_g$, $1>C_g$ such that  
\begin{equation}\label{condg1g2}
|{g^1}'(x)|\geq C_g|{g^2}'(x)|\geq C_g c_g \text{ for all }x\in\Omega 
\ \text{ and } \
\exists \, x_s\in\Omega: \ g^1(x_s)=g^2(x_s).
\end{equation} 
Note that under condition \eqref{condg1g2}, the range of $g^2$ is a subset of the range of $g^1$ and the set of points where the two functions coincide is a singleton. 
\begin{equation}\label{condg1g2_implications}
g^2(\Omega)\subseteq g^1(\Omega)
 \ \text{ and } \
\exists!\, x_s\in\Omega: \ g^1(x_s)=g^2(x_s)=:s_*.
\end{equation} 

This allows us to prove the following useful result.
\begin{lemma}\label{lem:Psi}
Under condition \eqref{condg1g2}, the mapping 
\[
\begin{aligned}
&\Psi: \{f_0\}+C_+^1([\ulu,\olu])\to C_*^1(\Omega), \quad
f\mapsto \log(f(g^1(x)))-\log(f(g^2(x)))
\\
&C_+^1([\ulu,\olu]):=\{f\in  C^1([\ulu,\olu])\, : \, f(s_0)=0, \ f+f_0\geq\underline{f}>0\}, \quad
\|f\|_{C_+^1[\ulu,\olu]}:= \|f'\|_{C[\ulu,\olu]}\\
&C_*^1(\Omega):=\{b\in C^1(\Omega)\, : \, b(x_s)=0\}, \quad
\|b\|_{C_*^1(\Omega)}:= \|b'\|_{C(\Omega)}
\end{aligned}
\]
is a  homeomorphism (that is, continuous and continuously invertible).
Moreover, the Lipschitz estimate
\begin{equation}\label{Psiinv}
\|\Psi^{-1}(\tilde{r}^1_I,\tilde{r}^2_I)-\Psi^{-1}(\tilde{r}^1_{II},\tilde{r}^2_{II})\|_{C_+^1([\ulu,\olu])}\leq 
C_\Psi \|(\tilde{r}^1_I,\tilde{r}^2_I)-(\tilde{r}^1_{II},\tilde{r}^2_{II})\|_{C_*^1(\Omega)}
\end{equation}
holds.
\end{lemma}

\begin{proof}
see the appendix.
\end{proof}
This can likely be extended to higher dimensions under more involved conditions on $g^1$, $g^2$.

\subsubsection*{case (c) mixed time-trace -- final time data, $\Gamma=\{x_0\}$}
Projecting the PDE onto the observation manifold we obtain (skipping the superscript ${}^1$)
\[
\begin{aligned}
&-\partial_t^\alpha h(t)+\triangle u(x_0,t) +r(x_0,t) = q(x_0) f(h(t))\quad && t\in(0,T)\\
&-\partial_t^\alpha u(x,T)+\triangle g(x)+r(x,T) = q(x) f(g(x))&& x\in\Omega.
\end{aligned}
\]
We assume $q(x_0)=q_1$ to be known, e.g., by using \eqref{qbar} with $\bar{x}=x_0$
and moreover make the assumption 
\begin{equation}\label{rangecondition}
u(\Omega\times(0,T))\subseteq h((0,T))=:[\ulu,\olu].
\end{equation}
From this, a fixed point scheme results as follows:
Given an iterate $(q^{(k)},f^{(k)})$, compute $f^{(k+1)}$ from
\begin{equation}\label{f_kp1}
q_1 f^{(k+1)}(h(t)) = -h'(t)+\partial_t^\alpha u^{(k)}(x_0,t) + q^{(k)} f^{(k)}(u^{(k)}_t(x_0,t)) \quad t\in(0,T)
\end{equation}
(where we have used the PDE to express $\triangle u$ via $\partial_t^\alpha u$ thus needing lower order  numerical differentiation) and then compute $q^{(k+1)}$ from
\begin{equation}\label{q_kp1}
\begin{aligned}
q^{(k+1)}(x)&=\Bigl(-\partial_t^\alpha u^{(k)}(x,T)+\triangle g(x)+r(x,T)\Bigr)/f^{(k+1)}(g(x)) \text{ or }\\
q^{(k+1)}(x)&=\Bigl(-\partial_t^\alpha \check{u}^{(k)}(x,T)+\triangle g(x)+r(x,T)\Bigr)/f^{(k+1)}(g(x)).
\end{aligned}
\end{equation}
Here $u^{(k)}$ solves \eqref{PDE} with $q=q^{(k)}$, $f=f^{(k)}$ and $\check{u}^{(k)}$ solves \eqref{PDE} with $q=q^{(k)}$, $f=f^{(k+1)}$.

\section{Reconstructions}

In this section we illustrate the performance of the devised reconstruction schemes in a spatially one-dimensional domain $\Omega=(0,L)$ over a time interval $(0,T)$ with $L=T=1$. For the numerical solution of the forward problem, we use a time stepping method based on the Abel integral for the Djrbashian-Caputo derivative, combined with a finite difference discretization in space. Due to the $\alpha$ dependent convergence rate, an optimal (temporal) mesh width would have to depend on $\alpha$. We here simply increased spatial and temporal resolution until the discetization error was not visible in the reconstruction results.

We first of all show reconstruction results in the noiseless case with the fixed point schemes in the final time (b) (iterates 2, 3, 9) and mixed (c) (iterates 3, 10, 39) observation settings.

\newpage
\includegraphics[width=\textwidth]{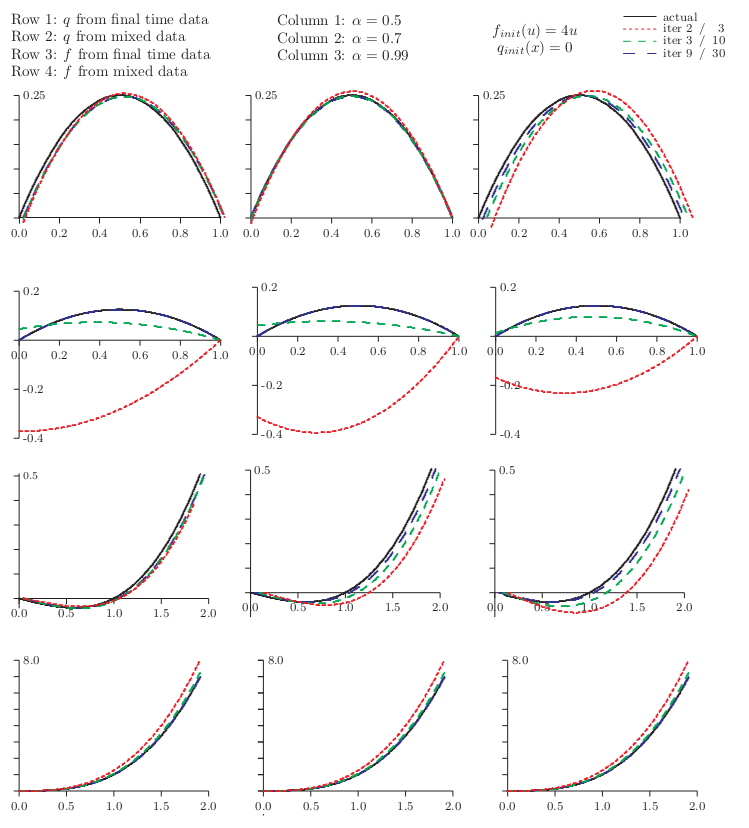}

\bigskip

The influence of $\alpha$ on the performance of the reconstruction schemes is not obvious. While a larger $\alpha$ should lead to a faster decay of $u$ and therefore to a smaller contraction constant, the interplay between $q$ and $f$ appears to be much more complex than the upper bounds in our analysis of Section~\ref{sec:uniquenessb} (and its extrapolation to the fractional case, where we expect a slower decay for smaller $\alpha$, thus a worse contraction constant) predicts.  

A comparison with the frozen Newton schemes shows that while the expected convergence rate is the same as for the fixed point schemes (namely linear) they are de facto by far slower.
The $L^2$ error in the (b) final time observation case was\hfill\break\\ 
\qquad 0.010352 for $q$ / 0.000766 for $f$ with $\alpha=0.5$ and \hfill\break\\
\qquad 0.012467 for $q$ / 0.000837 for $f$ with $\alpha=0.7$ after 100 frozen Newton iterations.\hfill\break\\
The fixed point scheme reached this level of precision already after three steps; at this point the $L^2$ error was \hfill\break\\ 
\qquad 0.000207 for $q$ / 0.000572 for $f$ with $\alpha=0.5$ and\hfill\break\\
\qquad 0.000151 for $q$ / 0.000466 for $f$ with $\alpha=0.7$.\hfill\break\\
It went down to \hfill\break\\
\qquad 0.000200 for $q$ / 0.000162 for $f$ with $\alpha=0.5$ and \hfill\break\\
\qquad 0.000191 for $q$ / 0.000133 for $f$ with $\alpha=0.7$ within 8 fixed point steps.\hfill\break\\
The corresponding numbers for the (c) mixed observation case were\hfill\break\\
\qquad 0.072306 for $q$ / 0.000822 for $f$ with $\alpha=0.5$ and\hfill\break\\
\qquad 0.072388 for $q$ / 0.000881 for $f$ with $\alpha=0.7$ after 100 frozen Newton steps and\hfill\break\\
\qquad 0.010057 for $q$ / 0.000904 for $f$ with $\alpha=0.5$ and\hfill\break\\
\qquad 0.011297 for $q$ / 0.001024 for $f$ with $\alpha=0.7$ after 25 fixed point steps.\hfill\break\\
For $\alpha=0.9$ and larger, frozen Newton failed to converge from the given starting values of $q$ and $f$ due to its obviously smaller convergence radius.

We finally report on reconstructions from noisy data.
Among the problem settings considered here, the least robust one to noise was time trace observations (a), which is intuitive since orthogonality of the observation manifold to variability of $q(x)$ makes its construction severely ill-posed, while reconstruction of $f$ should be neutral to whether the observation is done in space or in time direction.
Even after reducing the size $L$ of the spatial domain from unity to $0.2$, time trace Newton for (a)  only allowed 0.1 per cent noise (note that no fixed point scheme is available for this case). Larger noise quickly led to a failure of the solver, which is obviously sensitive to changes in the nonlinearity $f$.
The final time and mixed observation cases (b) and (c) were more robust: Both Newton and fixed point scheme allowed 1 and more per cent noise. 
In Table~\ref{tab:fiti-fxp_noise} we show $L^2$ errors with several noise levels and several values of the fractional differentiation order. Convergence as the noise level tends to zero is clearly visible in each of the cases. However, as in the noiseless reconstructions there is no monotonicity with respect to $\alpha$.  
\begin{table}
\begin{tabular}[h]{|l||c|c|c|c|}
\hline
			  & $\alpha=0.5$ & $\alpha=0.7$ & $\alpha=0.9$  & $\alpha=0.99$\\
\hline \hline 
$\delta=0.1\%$&0.000494		&0.002114		&0.000671		& 0.002378\\
$\delta=0.3\%$&0.005638		&0.006119		&0.004462		& 0.006115\\
$\delta=1\%$  &0.012843		&0.051335		&0.020535		& 0.015470\\
\hline \hline 
$\delta=0.1\%$&0.001089		&0.000888		&0.000731		& 0.001716\\
$\delta=0.3\%$&0.004603		&0.002750 		&0.003954		& 0.002461\\
$\delta=1\%$  &0.016325		&0.034767		&0.016931		& 0.011021\\
\hline
\end{tabular}
\caption{$L^2$ reconstruction errors in $q$ (top) and $f$ (bottom) with noisy data for different differentiation orders $\alpha$ for fixed point scheme \eqref{fixedpoint_finaltime_overall} with final time data (b).
\label{tab:fiti-fxp_noise}}
\end{table}

\section{Uniqueness and convergence 
of fixed point scheme 
in final time case (b) with $\alpha=1$}\label{sec:uniquenessb}

Consider the fixed point operator underlying the iterative scheme \eqref{fixedpoint_finaltime}, that is, 
\begin{equation}\label{T}
\mathbb{T}:(q,f)\mapsto 
\left(\begin{array}{c}
\Psi^{-1}\bigl(\log|\mathbb{S}_1(q,f)|,\log|\mathbb{S}_2(q,f)|\bigr)\\
\mathbb{S}_1(q,f)/\Psi^{-1}\bigl(\log|\mathbb{S}_1(q,f)|,\log|\mathbb{S}_2(q,f)|\bigr)
\end{array}\right),
\end{equation}
cf. Lemma~\ref{lem:Psi}, where 
\begin{equation}\label{S}
\mathbb{S}:(q,f)\mapsto(u^i_t(T)+\triangle g^i-r^i)_{i\in\{1,2\}} \text{ where $u^i$ solves \eqref{PDE} with $r=r^i$}.
\end{equation}
It follows from Lemma~\ref{lem:Psi} that $\Psi^{-1}$ is Lipschitz continuous, cf. \eqref{Psiinv}
and so is $\log(|\cdot|:\mathcal{B}_{\tilde{\rho}}^{C^1(\Omega)}(q_{\textrm{act}}f_{\textrm{act}}(u^i_{\textrm{act}}(T))\to C^1(\Omega)$ for $\tilde{\rho}>0$ small enough (noting that $\mathbb{S}_i(q,f)=u^i_t(T)+\triangle g^i-r^i$ is close to ${u_{\textrm{act}}}^i_t(T)+\triangle g^i-r^i=q_{\textrm{act}}f_{\textrm{act}}(u^i_{\textrm{act}}(T))\geq\underline{q}\underline{f}>0$ for $(q,f)$ close to $(q_{\textrm{act}},f_{\textrm{act}})$, see Theorem~\ref{thm:contractivityS}.

Contractivity of $\mathbb{T}$ and therewith uniqueness according to Banach's Contraction Principle can therefore be concluded (cf. Corollary~\ref{contractivityT}) from contractivity of $\mathbb{S}$ for sufficiently large $T$, which we state and prove now.

\begin{theorem}\label{thm:contractivityS}
Let $\alpha=1$, $\Omega\subseteq\mathbb{R}^1$ be a union of finitely many intervals, $\triangle u_0\in H^s(\Omega)$ for some $s>\frac12$, $r^i\equiv0$, $i\in\{1,2\}$ and let $\mu=\lambda_1/2$, where $\lambda_1$ is the smallest eigenvalue of the Dirichlet Laplacian on $\Omega$.\hfill\break\\
Then for any $\underline{q}\geq0$, $\rho_0$,  $\rho_1$, $\rho_2>0$ there exist  constants $C>0$ (large enough), $\rho_{qf},\,\rho_q>0$ (small enough) independent of $T$ such that for the mapping $\mathbb{S}$
the estimate 
\begin{equation}\label{LipschitzS}
\|\mathbb{S}(q_I,f_I)-\mathbb{S}(q_{II},f_{II})\|_{W^{1,\infty}(\Omega)}^2
\leq C e^{-\mu T}(\|q_I-q_{II}\|_{L^\infty(\Omega)} + \|f_I-f_{II}\|_{W^{1,\infty}_{\text{aff}}(\mathbb{R})})
\end{equation}
holds for any $q_I,q_{II}\in\mathcal{B}_{\rho_0}^{H^s(\Omega)}(0)$, and any $f_I\in \mathcal{B}_{\rho_1}^{W^{1,\infty}(\mathbb{R})}(0)$, $f_{II}\in \mathcal{B}_{\rho_2}^{W^{1,\infty}_{\textrm{aff}}(\mathbb{R})}(0)$ such that 
$\|\frac{1}{T|\Omega|}\int_0^T\int_\Omega q_{II}(x)f_{II}'(u_{II}(x,t))\,dx\,dt -q_{II}f_{II}'(u_{II})\|_{L^\infty(0,T;L^1(\Omega))}<\rho_{qf}$ and 
$\|q_I-q_{II}\|_{L^\infty(\Omega)}\leq \rho_q$.
\end{theorem}
\begin{proof}
\underline{Step 1: Decay estimate for $u^i_{I\,t}(t)$, $u^i_{II\,t}(t)$:}\hfill\break\\
Skipping the subscripts $I$, $II$ for ease of notation, we define $v=\frac12 (u_t)^2$ for a solution $u$ of \eqref{PDE} with $r\equiv0$. It is readily checked that 
\[
\begin{aligned}
&v_t-\triangle v =u_t(u_t-\triangle u)_t-|\nabla u_t|^2 = - qf'(u)(u_t)^2 -|\nabla u_t|^2\leq0, \ t>0\\
&v(t=0)=v_0:=\bigl(\triangle u_0-qf(u_0)\bigr)^2
\end{aligned}
\]
and therefore, by the maximum principle, $0\leq v\leq\bar{v}$ where $\bar{v}$ solves
\[
\bar{v}_t-\triangle \bar{v} = 0, \ t>0 \qquad \bar{v}(t=0)=v_0.
\]
From this we can conclude the following decay estimates on $u_t$.
\begin{equation}\label{decay_u}
\begin{aligned}
\|u_t(t)\|_{L^4(\Omega)}^2 &= 2 \|v(t)\|_{L^2(\Omega)} \leq 2 \|\bar{v}(t)\|_{L^2(\Omega)}
=2 \left(\sum_{j=1}^\infty e^{-2\lambda_j t}\langle v_0,\varphi_j\rangle^2\right)^{1/2}\\
&\leq 2 e^{-\lambda_1 t} \|v_0\|_{L^2(\Omega)}\\[1ex]
\|u_t(t)\|_{L^\infty(\Omega)}^2 &= 2 \|v(t)\|_{L^\infty(\Omega)} \leq 2 \|\bar{v}(t)\|_{L^\infty(\Omega)}
\leq 2 C_{H^{\bar{s}},L^\infty}^\Omega \left(\sum_{j=1}^\infty e^{-2\lambda_j t}\lambda_j^{2\bar{s}} \langle v_0,\varphi_j\rangle^2\right)^{1/2}
\\&\leq 2 C_{H^{\bar{s}},L^\infty}^\Omega e^{-\lambda_1 t} \|v_0\|_{H^{\bar{s}}(\Omega)}
\end{aligned}
\end{equation}
with an eigensystem $(\lambda_j,\varphi_j^0)_{j\in\mathbb{N}}$ of $-\triangle_B$ and $\bar{s}=\min\{1,s\}>\frac12$ so that $H^{\bar{s}}(\Omega)$ continuously embeds into $L^\infty(\Omega)$. With this choice, $H^{\bar{s}}(\Omega)$ is also a Banach Algebra 
and therefore 
\begin{equation}\label{est_v0}
\|v_0\|_{H^{\bar{s}}(\Omega)}\leq \|\triangle u_0-qf(u_0)\|_{H^{\bar{s}}(\Omega)}^2
\end{equation}
Alternatively, we could consider the PDE for $\bar{u}(x,t):=e^{\mu t}u_t(t)$
\begin{equation}\label{PDE_ubar}
\bar{u}_t-\triangle \bar{u}+(qf'(u)-\mu)\bar{u}=0,
\ t>0, \quad  \ \bar{u}(0)=\triangle u_0-qf(u_0)
\end{equation}
and apply estimates analogous to those in step 3 below. This would have the advantage of generalising to a higher dimensional domain $\Omega$ and a more general elliptic operator $\mathbb{L}$ in place of the Laplacian, as long as a Green's function for $\partial_t-\mathbb{L}+(\bar{c}-\mu)$ exists and has finite norms as appearing in step 3.
However, we are anyway limited to the 1-d setting by Lemma~\ref{lem:Psi}.
\hfill\break\\[1ex]
\underline{Step 2: PDE for difference $z$ and $w:=z_t$:}\hfill\break\\
In view of the fact that $\mathbb{S}(q_I,f_I)-\mathbb{S}(q_{II},f_{II})=(u^i_{I\,t}(T)-u_{II\,t}(T))_{i\in\{1,2\}}$, 
we abbreviate  $z^i:= u^i_I-u^i_{II}$, $w^i:=u^i_{I\,t}-u_{II\,t}$ and obtain from \eqref{PDE}, skipping the superscript $i$ for simplicity, 
\[
\begin{aligned}
&z_t-\triangle z = q_{II}f_{II}(u_{II})-q_If_I(u_I), \ t>0\\
&z(t=0)=0
\end{aligned}
\]
and 
\[
\begin{aligned}
&w_t-\triangle w = q_{II}f_{II}'(u_{II})u_{II,t}-q_If_I'(u_I)u_{I,t}, \ t>0\\
&w(t=0)=w_0:=q_{II}f_{II}(u_0)-q_If_I(u_0)
\end{aligned}
\]
where
\[
\begin{aligned}
&q_{II}f_{II}(u_{II})-q_If_I(u_I)\\
&=
(q_{II}-q_I)f_{II}(u_{II})
+q_I\bigl(f_{II}(u_{II})-f_{II}(u_I)
	+(f_{II}-f_I)(u_{I})\bigr)\\
&=
(q_{II}-q_I)f_{II}(u_{II})
-q_I\int_0^1 f_{II}'(u_{II}+\theta z)\, d\theta\, z
+q_I(f_{II}-f_I)(u_{I})\bigr)
\\[1ex]
&q_{II}f_{II}'(u_{II})u_{II,t}-q_If_I'(u_I)u_{I,t}\\
&=(q_{II}-q_I)f_{II}'(u_{II})u_{II,t}
+q_I\bigl(-\int_0^1 f_{II}''(u_{II}+\theta z)\, d\theta\, z u_{II,t}
	-f_{II}'(u_{I})w
	+(f_{II}'-f_I')(u_{I})u_{I,t}\bigr)
\end{aligned}
\]
Thus, with 
\[
\begin{aligned}
&c^z:=q_I\int_0^1 f_{II}'(u_{II}+\theta z)\, d\theta\geq0, \quad
c^w:=q_I f_{II}'(u_{I})\geq0, \quad
m:=q_I\int_0^1 f_{II}''(u_{II}+\theta z)\, d\theta\, u_{II,t}, \\
&r^z:= (q_{II}-q_I)f_{II}(u_{II})
+q_I(f_{II}-f_I)(u_{I})\bigr)\\
&r^w:= (q_{II}-q_I)f_{II}'(u_{II})u_{II,t}
	+q_I(f_{II}'-f_I')(u_{I})u_{I,t},
\end{aligned}
\]
we have
\begin{equation}\label{PDE_z}
\begin{aligned}
&z_t-\triangle z+c^z z=r^z, \ &t>0, \quad &z(0)=0\\
&w_t-\triangle w+c^w w=r^w-mz, \ &t>0, \quad &w(0)=w_0
\end{aligned}
\end{equation}
with homogeneous Dirichlet boundary conditions,
where due to \eqref{estfval}, \eqref{decay_u}, $m$, $r^z$, $r^w$, $w_0$ satisfy the estimates
\begin{equation}\label{decay_mr}
\begin{aligned}
&\|m(t)\|_{L^\infty(\Omega)} \leq C_1 e^{-(\lambda_1/2) t},\\
&\|r^z(t)\|_{L^\infty(\Omega)} \leq C_2 
(\|q_{II}-q_I\|_{L^\infty(\Omega)}+\|f_{II}-f_I\|_{W^{1,\infty}_{\textrm{aff}}(\mathbb{R})},\\
&\|r^w(t)\|_{L^\infty(\Omega)} \leq C_3 e^{-(\lambda_1/2) t} 
(\|q_{II}-q_I\|_{L^\infty(\Omega)}+\|(f_{II}-f_I)'\|_{L^\infty(\mathbb{R})}),\\
&\|w_0\|_{L^\infty(\Omega)} \leq C_4 
(\|q_{II}-q_I\|_{L^\infty(\Omega)}+\|f_{II}-f_I\|_{W^{1,\infty}_{\textrm{aff}}(\mathbb{R})}),
\end{aligned}
\end{equation}
with $C_1$, $C_2$, $C_3$, $C_4$ depending only on $\rho_0$, $\rho_1$, $\rho_2$, cf. \eqref{est_v0}.
In view of \eqref{decay_mr} we multiply the second line of \eqref{PDE_z} by $e^{\mu t}$ with $\mu=\min\{\underline{c},\lambda_1/2\}$, to obtain for 
$\bar{w}(x,t):=e^{\mu t}w(x,t)$, 
$\bar{m}(x,t):=e^{\mu t}m(x,t)$, $\bar{r}(x,t):=e^{\mu t}r(x,t)$
\begin{equation}\label{PDE_wbar}
\bar{w}_t-\triangle \bar{w}+(c-\mu)\bar{w} -\bar{m}\, z
=\bar{r},
\ t>0, \quad  \ \bar{w}(0)=w_0.
\end{equation}
\underline{Step 3: Estimate on $\bar{w}(T)=e^{\mu T} (u^i_{I\,t}(T)-u_{II\,t}(T)$:}
With the Green's function $G$ of $\partial_t-\triangle + \bar{c}^w$, cf, e.g. \cite{Cannon:1984},
for $\bar{c}^w=\frac{1}{T|\Omega|}\int_0^T\int_\Omega c^w(x,t)\,dx\,dt$
we can write 
\[
\begin{aligned}
&\|\bar{w}(T)\|_{C^k(\Omega)} = 
\|\int_\Omega G(\cdot,y,T)w_0(y)\,dy+\int_0^T G(\cdot,y,T-\tau)\bar{r}_{tot}(y,\tau)\,dy\, d\tau\|_{C^k(\Omega)}\\
&\leq \|w_0\|_{L^1(\Omega)} \|G(\cdot,\cdot,T)\|_{L^\infty(\Omega;C^k(\Omega))}
+\|\bar{r}_{tot}\|_{L^\infty(0,T;L^1(\Omega))}
\|G\|_{L^1(0,T;L^\infty(\Omega;C^k(\Omega)))}
\end{aligned}
\]
where $\bar{r}_{tot}=\bar{r}^w-\bar{m}z+(\bar{c}^w-c^w)w$ and thus
\[
\begin{aligned}
\|\bar{r}_{tot}\|_{L^\infty(0,T;L^1(\Omega))}
\leq &
\|\bar{r}^w\|_{L^\infty(0,T;L^1(\Omega))}
+\|\bar{m}\|_{L^\infty(0,T;L^1(\Omega))} \|z\|_{L^\infty(0,T;L^\infty(\Omega)}
\\&+\|\bar{c}^w-c^w\|_{L^\infty(0,T;L^1(\Omega))} \|w\|_{L^\infty(0,T;L^\infty(\Omega)}.
\end{aligned}
\]
This holds for arbitrary $k\in\mathbb{N}_0$ and reflects the well-known infinite smoothing properties of the solution operator of the heat equation. For our purposes $k=1$ suffices.
Applying similar estimates with $T$ replaced by $t$ and $k=0$ we obtain
\[
\begin{aligned}
&\|\bar{z}\|_{L^\infty(0,t;L^\infty(\Omega))} 
\leq 
\|\bar{r}^v\|_{L^\infty(0,T;L^1(\Omega))}
\|G\|_{L^1(0,T;L^\infty(\Omega;L^\infty(\Omega)))}\\[1ex]
&\|\bar{w}\|_{L^\infty(0,t;L^\infty(\Omega))} 
\leq (\|w_0\|_{L^\infty(\Omega)} 
+\|\bar{r}_{tot}\|_{L^\infty(0,T;L^1(\Omega))})
\|G\|_{L^1(0,T;L^\infty(\Omega;L^\infty(\Omega)))},
\end{aligned}
\]
thus with 
\[
\begin{aligned}
&C_G:=\|G\|_{L^1(0,T;L^\infty(\Omega;L^\infty(\Omega)))}\\
&d(q,f):=\|\bar{r}^w\|_{L^\infty(0,T;L^1(\Omega))}
+\|\bar{m}\|_{L^\infty(0,T;L^1(\Omega))} 
\|\bar{r}^v\|_{L^\infty(0,T;L^1(\Omega))}C_G
\end{aligned}
\]
we have 
\[
\|\bar{w}\|_{L^\infty(0,t;L^\infty(\Omega))}\leq
\frac{C_G\,d(q,f)}{1-C_G\|\bar{c}^w-c^w\|_{L^\infty(0,T;L^1(\Omega))}}
\]
and, with \eqref{decay_mr} and some constant $C_5$ independent of $T$,
\[
\begin{aligned}
\|\bar{w}(T)\|_{C^k(\Omega)} 
&\leq \|w_0\|_{L^1(\Omega)} \|G(\cdot,\cdot,T)\|_{L^\infty(\Omega;C^k(\Omega))}
+\frac{1}{1-C_G\|\bar{c}^w-c^w\|_{L^\infty(0,T;L^1(\Omega))}}d(q,f)\\
&\leq
\frac{C_5}{1-C_G\|\bar{c}^w-c^w\|_{L^\infty(0,T;L^1(\Omega))}} 
(\|q_{II}-q_I\|_{L^\infty(\Omega)}+\|f_{II}-f_I\|_{W^{1,\infty}(\mathbb{R})}).
\end{aligned}
\]
Note that smallness of $\|\bar{c}^w-c^w\|_{L^\infty(0,T;L^1(\Omega))}$ follows from the assumed smallness of \hfill\break\\
$\|\frac{1}{T|\Omega|}\int_0^T\int_\Omega q_{II}(x)f_{II}'(u_{II}(x,t))\,dx\,dt -q_{II}f_{II}'(u_{II})\|_{L^\infty(0,T;L^1(\Omega))}$ for $q_I$ sufficiently close (wrt. the $L^\infty(\Omega)$ norm) to $q_{II}$.
\end{proof}

\medskip

Contractivity of $\mathbb{T}$ therefore follows from contractivity of $\mathbb{S}$ for sufficiently large $T$, as stated in  Theorem~\ref{thm:contractivityS}.

\begin{corollary}\label{contractivityT}
Under the conditions of Lemma~\ref{lem:Psi} and Theorem~\ref{thm:contractivityS} with $(q_{II},f_{II}):=(q_{\textrm{act}},f_{\textrm{act}})$, there exist $\rho>0$ (sufficiently small), $T>0$ (sufficiently large) and $c\in(0,1)$ such that the fixed point operator $\mathbb{T}$ satisfies the contractivity estimate 
\[
\| \mathbb{T}(q,f)- \mathbb{T}(q_{\textrm{act}},f_{\textrm{act}})\|_{L^\infty(\Omega),W^{1,\infty}([\ulu,\olu])}
\leq c
\|(q,f)- (q_{\textrm{act}},f_{\textrm{act}})\|_{L^\infty(\Omega),W^{1,\infty}([\ulu,\olu])}
\]
in $\mathcal{B}_\rho^{L^\infty(\Omega),W^{1,\infty}([\ulu,\olu])}(q_{\textrm{act}},f_{\textrm{act}})$. Therefore $(q,f)$ is locally \footnote{in $\mathcal{B}_\rho^{L^\infty(\Omega),W^{1,\infty}([\ulu,\olu])}(q_{\textrm{act}},f_{\textrm{act}})$} uniquely determined by the final time measurements $(b) \ g^i(x) = u^i(x,T), \ x\in \Omega$, $i\in\{1,2\}$.
\end{corollary}

As a consequence, we obtain convergence of the iteration scheme \eqref{fixedpoint_finaltime_overall} with exact data.

\medskip

\subsubsection{Convergence with noisy data}
To handle noisy data $\tilde{y}$ in place of $y$, we combine \eqref{fixedpoint_finaltime} with a pre-smoothing step that constructs $y^\delta\in W^{2,\infty}(\Omega)\subseteq C^1(\Omega)$ such that 
\begin{equation}\label{delta}
\|y^\delta-y\|_{W^{2,\infty}(\Omega)}\leq\delta,
\end{equation}
and use $y^\delta$ in place of $y$ in the fixed point scheme.
This means that  we apply the Picard iteration 
\begin{equation}\label{PicardTdelta}
(q^{(k+1)\,\delta},f^{(k+1)\,\delta})=\mathbb{T}^\delta(q^{(k)\,\delta},f^{(k)\,\delta}), \quad (q^{(0)\,\delta},f^{(0)\,\delta}) =(q^{(0)},f^{(0)})
\end{equation}
with $\mathbb{T}^\delta$ defined by $\mathbb{T}^\delta(q,f):=(\Psi^\delta)^{-1}\bigl(\log(|\mathbb{S}_1^\delta(q,f)|,\log(|\mathbb{S}_2^\delta(q,f)|\bigr)$.

To see that \eqref{Psiinv} remains valid with a Lipschitz constant $C_\Psi$ that can be chosen independently of $\delta$, we inspect the assumptions of Lemma~\ref{lem:Psi} with respect to their dependence on the data. 
If condition \eqref{condg1g2} is satisfied for the exact data $(g^1,g^2)=F(q_{\textrm{act}},f_{\textrm{act}})$ then it remains valid for $(g^1,g^2)=y^\delta$ with constants $C_g,\,c_g>0$ that can be chosen independent of $\delta$ and a unique intersection point $x_s^\delta$ that satisfies $|x_s^\delta-x_s|\leq C_s\delta$ with $x_s$ being the exact intersection point and $C_s$ independent of $\delta$.
As a consequence, also \eqref{condg1g2_implications} remains valid.

Dependence of $\mathbb{S}$ on the data is much simpler, occurring only through the $\triangle {g^\delta}^i$ term and therefore cancelling out in the proof of contractivity Theorem~\ref{thm:contractivityS}, which thus remains valid without any changes.

Therefore, we can invoke Corollary~\ref{contractivityT}, which implies that the sequence $(q^{(k)\,\delta},f^{(k)\,\delta})_{k\in\mathbb{N}}$ defined by \eqref{PicardTdelta} converges linearly to some fixed point $(q^\delta,f^\delta)\in X$ of $\mathbb{T}^\delta$. 
The estimate
\[\begin{aligned}
\|(q^\delta,f^\delta)-(q_{\textrm{act}},f_{\textrm{act}})\|_X
=&\|\mathbb{T}^\delta(q^\delta,f^\delta)-\mathbb{T}(q_{\textrm{act}},f_{\textrm{act}})\|_X\\
\leq&\|\mathbb{T}^\delta(q^\delta,f^\delta)-\mathbb{T}^\delta(q_{\textrm{act}},f_{\textrm{act}})\|_X
+\|(\mathbb{T}^\delta-\mathbb{T})(q_{\textrm{act}},f_{\textrm{act}})\|_X,
\end{aligned}\]
where due to the reasoning above $\|(\mathbb{T}^\delta-\mathbb{T})(q_{\textrm{act}},f_{\textrm{act}})\|_X\leq C\delta$, together with contractivity (with a $\delta$ independent constant $c\in(0,1)$) of $\mathbb{T}^\delta$ yields the stability bound
\[
\|(q^\delta,f^\delta)-(q_{\textrm{act}},f_{\textrm{act}})\|_X\leq \frac{C}{1-c} \delta.
\]
Note that the iteration \eqref{PicardTdelta} does not require any regularization nor early stopping; it suffices to counteract ill-posedness by smoothing the data.


\subsection*{Acknowledgment}
\vskip-4pt
\noindent
The work of the first author was funded in part by the Austrian Science Fund (FWF) 
[10.55776/P36318]. 

\noindent
The work of the second author was supported in part by the
National Science Foundation through award {\sc dms}-2111020.


\input{appendix_proofLemPsi}

\input{appendix_discr} 

\end{document}

%% file: intro.tex
\section{Introduction}

In this paper, we study reconstruction of both $q$ and $f$ 
in a reaction-subdiffusion equation with the specific form
\begin{equation}\label{eqn:qf_parabolic}
\frac{\partial u}{\partial t} - \triangle u + q(x)f(u) = 0
\end{equation}
defined in a cylindrical domain $\Omega\times(0,T)$ where $\Omega$
is a closed bounded domain in $\mathbb{R}^d$ with smooth boundary
$\partial\Omega$.
Here $q(x)$ is a spatially-dependent coefficient and the reaction term $f(u)$
depends only on the solution $u$.

Such equations are ubiquitious in applications throughout the
biological and physical sciences.
For example, in population dynamics the famous Verhulst
ordinary differential equation model sets $f(u) = au -b u^2$ to
separate out the birth and death processes.
This was extended to spatial regions and to include diffusion by Skellum
in 1951, \cite{Skellum:1951} and was one of the first models to incorporate
random walk processes into the subject.
The model therefore parallels the Fisher model for an advantageous gene
in a population introduced slightly earlier \cite{Fisher:1937,Fisher:book}.
Indeed, the general model \eqref{eqn:qf_parabolic} with $q$ constant and $f$ of the form
$f(u) = a(u -u^2)$ is now commonly referred to as the Fisher equation.
There is a vast literature in this application context for the equation
and allowing a spatial dependence on the coefficients is clearly
of importance.
The generalisation to include diffusion and hence spatial dependence,
adds considerabe mathematical complexity including the possibility of
travelling wave solutions of the above equation.
See, for example, \cite{Murray,AronsonWeinberger}.
The form of the right hand side term here balances growth and decay with
different power-law terms and clearly there are many possible modifications
from a modelling standpoint.
These might include different rates by using exponents other than linear and
quadratic terms as well as higher order and fractional powers involving $u$,
or indeed a functional form that goes beyond just power law behaviour.

In perhaps the other most well-known application, that of chemical reactions,
a similar form is used, \cite{Grindrod:1996, Murray:2002}.
A cubic polynomial for $f$ in terms of $u$ appears in the Rayleigh-B\'ernard
model of convection and in the FitzHugh--Nagumo model of 1961
modelling impulses along a nerve, \cite{FitzHugh:1961}
and also in the Zeldovich--Frank--Kamenetsky--Equation.
More general models of the ZFK type involve
exponential functions: for example, $f(u) = \alpha\,u(1-u)e^{-\beta (1-u)}$,
\cite{FrankZeldovich}.
Thus there is an application-driven reason to go beyond merely seeking
$f(u)$ as a combination of powers of $u$.

Thus, although these models are based on specific polynomial
assumptions there have been numerous modifications including
assuming more general polynomials and with possibly non-integer powers.
Our paradigm here allows us to recover the  form of the reaction term
from measured data; and, in fact, in addition to allow for
a coupled inhomogeneous spatial environment that itself has to recovered.
We note that this type of coupling is a major factor in the difficulty
of determining both $q$ and $f$.

\medskip
Undetermined coefficient  problems for parabolic equations have a history
going back more than 40 years.
For the linear equation $u_t - \triangle u + q(x)u = 0$, 
Pierce showed a uniqueness result for recovering $q(x)$
in one spatial variable from time trace data $u(x_0,t)$
with $x_0$ on the boundary \cite{Pierce:1979}. He also indicated a reconstruction method.
For the case of final time data $u(x,T)$, Rundell proved unique identifiability
for $q$ and showed it could be determined by an iterative process
from a fixed point scheme, \cite{Rundell:1983,Rundell:1987}.
For the case where $q(x)=1$ and the nonlinear term $f$ was unknown
Pilant and Rundell in a series of papers 
\cite{PilantRundell:1986,PilantRundell:1987,PilantRundell:1988},
showed uniqueness and  gave reconstruction methods for $f(u)$. 
See also the book by Isakov, \cite{Isakov:2017}.

More recently, due to its relevance in modeling, time-fractional 
versions of these problems have been considered.
Also in this paper we will include the subdiffusion case where our model 
extends the differential equation for $0<\alpha\leq 1$ by
\begin{equation}\label{PDE_intro}
\partial_t^\alpha u-\Delta u +q(x) f(u) = r 
\end{equation}
The literature in this direction and from an inverse problems standpoint
is growing rapidly and we reference the book \cite{BBB} as a starting point.

We here only give a small selection of related examples.
The most closely related ones are identification of the spatially dependent potential $q(x)$ in \eqref{PDE_intro} for given $f(u)$ in \cite{JinRundell:2012b,KaltenbacherRundell:2019b} or the determination of $f(u)$ in \eqref{PDE_intro} for given $q(x)$ in \cite{KaltenbacherRundell:2019b}.
Uniqueness and reconstructibility results for the pair $\{a(x),f(u)\}$ in
 $u_t^\alpha - \bigl(a(x)u_x)_x - f(u)$
can be found in \cite{KaltenbacherRundell:2020c}. 

From the above it can be seen that our interest here lies in the
inverse problem of recovering {\it both\/} the spatial term $q(x)$
{\it and\/} the specific form of the nonlinearity $f(u)$.
This is clearly an essential problem as while the general form of the model
might be prescribed in an application the individual components are likely
to be unknown and must be determined from additional (overposed) data
measurements.

The additional difficulty over previous combinations noted above
is the fact that the two unknowns are directly coupled as a product
as distinct from an unknown diffusion coefficient $a(x)$ and reaction $f(u)$.

Also the type of overposed data is highly relevant for reconstructibility of $q(x)$ and $f(u)$.
We here study on one hand measurements inside $\Omega$ at a fixed time instance $T$, on the other hand, boundary observations over a time interval $(0,T)$, as well as a combination thereof.

The plan of the paper is as follows.
In section~\ref{sec:intro} we formulate the inverse problem in the three above mentioned observation settings.
Section~\ref{sec:analysis} is devoted to an analysis of the forward problem and a proof of well-definedness and Lipschitz continuous Fr\'{e}chet differentiability of the operator mapping $(q,f)$ to the solution of \eqref{PDE_intro} with initial and boundary conditions.
In section~\ref{sec:iterations} we develop some iterative reconstruction methods for each of these scenarios. Some of them result from  fixed point schemes that are based on a projection of the PDE on the observation manifold; for comparison, we also consider Newton type schemes. In the final time data case with $\alpha=1$ and $\Omega\subseteq\mathbb{R}^1$ we provide a convergence result of the fixed point iteration and prove uniqeness in section~\ref{sec:uniquenessb}.

%% file: appendix_proofLemPsi.tex
\section*{Appendix A: Proof of Lemma~\ref{lem:Psi}}
\begin{proof}
Having in mind $\tilde{f}=\log(|f|)-\log(f_0)$, we define the linear operator 
\[
\begin{aligned}
&\Psi_{lin}:C_0^1([\ulu,\olu])\to C_*^1(\Omega), \quad 
\tilde{f}\mapsto \tilde{f}\circ g^1 - \tilde{f}\circ g^2
\\
&C_0^1([\ulu,\olu]):=\{\tilde{f}\in  C^1([\ulu,\olu])\, : \, \tilde{f}(s_0)=0\}, \quad
\|\tilde{f}\|_{C_0^1[\ulu,\olu]}:= \|\tilde{f}'\|_{C[\ulu,\olu]}.
\end{aligned}
\]
By differentiating  the above we obtain, choosing 
$x_*\in \text{argmax}|\tilde{f}'\circ g^1|$, and using \eqref{condg1g2},
\begin{equation}\label{Psilin_injective}
\begin{aligned}
\|(\Psi_{lin}\tilde{f})'\|_{C(\Omega)}
&\geq |{g^1}'(x_*)|\,|\tilde{f}'(g^1(x_*))|-|{g^2}'(x_*)|\,|\tilde{f}'(g^2(x_*))|\\
&= |{g^1}'(x_*)|\, \|\tilde{f}\|_{C_0^1[\ulu,\olu]}-|{g^2}'(x_*)|\,|\tilde{f}(g^2(x_*))|\\
&\geq (C_g-1)c_g \|\tilde{f}\|_{C_0^1[\ulu,\olu]}.
\end{aligned}
\end{equation}
This implies injectivity of $\Psi_{lin}$. To prove that it is also surjective, we fix $b\in C^1_*(\Omega)$ and construct a preimage under $\Psi_{lin}$ by the fixed point iteration
\[
\begin{aligned}
&{\tilde{f}'}_{(n+1)}(s)=\frac{b((g^1)^{-1}(s))+{g^2}'((g^1)^{-1}(s)){\tilde{f}'}_{(n)}(g^2(s))}{{g^1}'((g^1)^{-1}(s))}\\
&\tilde{f}_{(n+1)}(s_0)=0
\end{aligned}
\]
The corresponding fixed point operator can be easily shown to be contractive, due to the estimate (note the definition of the norm on $C_0^1[\ulu,\olu]$)
\[
\begin{aligned}
&\left\|\left(\frac{b+{g^2}'{\tilde{f}_I'}\circ g^2}{{g^1}'}\right)\circ(g^1)^{-1}
-\left(\frac{b+{g^2}'{\tilde{f}_{II}'}\circ g^2}{{g^1}'}\right)\circ(g^1)^{-1}\right\|_{C[\ulu,\olu]}\\
&=\left\|\frac{b+{g^2}'{\tilde{f}_I'}\circ g^2}{{g^1}'}
-\frac{b+{g^2}'{\tilde{f}_{II}'}\circ g^2}{{g^1}'}\right\|_{C(\Omega)}
=\left\|\frac{{g^2}'}{{g^1}'}(\tilde{f}_I'-\tilde{f}_{II}')\circ g^2\right\|_{C(\Omega)}\\
&\leq\left\|\frac{{g^2}'}{{g^1}'}\right\|_{C(\Omega)} \|\tilde{f}_I'-\tilde{f}_{II}'\|_{C[\ulu,\olu]}
\leq \frac{1}{C_g} \|\tilde{f}_I-\tilde{f}_{II}\|_{C_0^1[\ulu,\olu]}.
\end{aligned}
\]
Thus ${\tilde{f}'}_{(n)}$ converges to some $\tilde{f}$ in the Banach space $C_0^1[\ulu,\olu]$ and it is readily checked that $\Psi_{lin}\tilde{f}$. 
We have therefore verified surjectivity of $\Psi_{lin}\tilde{f}$. 
With this and the above estimate \eqref{Psilin_injective}, we have bijectivity and bounded invertibility of $\Psi_{lin}$.

Moreover, the nonlinear map 
\[
\Psi_{nl}: \{f_0\}+C_+^1([\ulu,\olu])\to C_0^1([\ulu,\olu]), \quad
f\mapsto \log(f)-\log(f_0)
\]
is Lipschitz continuous and Lipschitz continuously invertible with $\Psi_{nl}^{-1}(\tilde{f})=f_0\exp(\tilde{f})$, due to the estimates
\[
\begin{aligned}
&\|\Psi_{nl}(f_1)-\Psi_{nl}(f_2)\|_{C_+^1([\ulu,\olu])}\\
&=\|\tfrac{f_1'}{f_1}-\tfrac{f_2'}{f_2}\|_{C_+^1([\ulu,\olu])}
=\|\tfrac{1}{f_1}(f_1'-f_2')-\tfrac{f_2'}{f_1f_2}\int_{s_0}^\cdot(f_1'-f_2')\, ds\|_{C_+^1([\ulu,\olu])}\\
&\leq (\tfrac{1}{\underline{f}}+\tfrac{\olu-\ulu}{\underline{f}^2})\|f_1-f_2\|_{C_0^1([\ulu,\olu])}
\\[1ex]
&\|\Psi_{nl}^{-1}(\tilde{f}_1)-\Psi_{nl}^{-1}(\tilde{f}_2)\|_{C_0^1([\ulu,\olu])}\\
&=\|\exp(\tilde{f}_1)\, (\tilde{f}_1'-\tilde{f}_2') 
+ \tilde{f}_2'\int_0^1 \exp(\tilde{f}_2+\vartheta(\tilde{f}_1-\tilde{f}_2))\, d\vartheta\,
\int_{s_0}^\cdot (\tilde{f}_1'-\tilde{f}_2')\, ds\|_{C([\ulu,\olu])}\\
&\leq (\|\exp(\tilde{f}_1)\|_{C([\ulu,\olu])}+
(\olu-\ulu)\|\tilde{f}_2'\int_0^1 \exp(\tilde{f}_2+\vartheta(\tilde{f}_1-\tilde{f}_2))\, d\vartheta \|_{C([\ulu,\olu])})
\|\tilde{f}_1-\tilde{f}_2\|_{C_+^1([\ulu,\olu])}
\end{aligned}
\]
The claim of the lemma follows from the fact that $\Psi=\Psi_{nl}\circ \Psi_{lin}$.
\end{proof}

%% file: appendix_discr.tex
\section*{Appendix B: Discretization of iterative schemes}
\subsubsection{Frozen Newton}\label{sec:Newton_app}
To discretize computation of the Newton step $(\dq,\df)$ from $F'(q^{(0)},f^{(0)})(\dq,\df)=y-F(q^{(k)},f^{(k)})$ where $q^{(0)}\equiv c$, $f^{(0)}(u)=u$, we make an ansatz
\[
\dq(x)=\sum_{m=1}^M a_m\, \chi_m(x), \qquad
\df(u)=\sum_{n=1}^N b_n\, \theta_n(u)
\]
with basis functions $\chi_m$, $\theta_n$ and set up the semi-discretized Jacobian
\[
J_m=C\du^q_m, \quad m\in\{1,\ldots M\}, \qquad J_{M+n}=C\du^f_n, \quad n\in\{1,\ldots N\},
\]
where $\du^q_m$ solves \eqref{dPDE} with $\dq=\chi_m$, $\df=0$, and $\du^f_n$ solves \eqref{dPDE} with $\dq=0$, $\df=\theta_n$. 
The coefficients $\underline{a}=(a_1,\ldots,a_M)$, $\underline{b}=(b_1,\ldots,b_M)$ are then obtained by solving the linear system
\[
J\left(\begin{array}{c}\underline{a}\\ \underline{b}\end{array}\right) = y-F(q^{(k)},f^{(k)})
\]
e.g., by some regularised least squares methods.

To take into account the fixed value of $q$ at $\bar{x}$ according to \eqref{qbar}, we use basis functions $\chi_m(x)$ that vanish at $\bar{x}$ and an inital guess $q^{(0)}$ that satisfies $q^{(0)}(\bar{x})=\bar{q}$.

\subsubsection*{Fixed point scheme in case (b)}
In order to fix the values 
$q(\bar{x})=\bar{q}$, $f(0)=f_0$,
we consider the decompositions $\tilde{q}=\log(\bar{q})+\tilde{q}_\sim$, $\tilde{f}=\log(f_0)+\tilde{f}_\sim$,  and 
represent 
$\tilde{q}_\sim$ and $\tilde{f}_\sim$ by basis functions that vanish at $x=\bar{x}$ and $u=0$, respectively.\\

To numerically implement the fixed point schemes \eqref{fixedpoint_finaltime} or \eqref{fixedpoint_finaltime_overall}, we discretize 
\begin{equation}\label{basis_fixedpoint_finaltime}
\tilde{q}^{(k+1)}(x)=\sum_{m=1}^M a_M\, \tilde{\chi}_m(x), \qquad
\tilde{f}^{(k+1)}(u)=\sum_{n=1}^N b_n\, \tilde{\theta}_n(u)
\end{equation}
with basis functions $\tilde{\chi}_m$, $\tilde{\theta}_n$ and solve \eqref{fixedpoint_finaltime} or \eqref{fixedpoint_finaltime_overall} in a least squares sense.  

In case of \eqref{fixedpoint_finaltime}, this leads to the linear systems 
$A^f \vec{b}= \mathfrak{r}^f$, $A^q \vec{a}= \mathfrak{r}^q$  with 
\[
\begin{aligned}
&A^f_{jn} =\int_\Omega \Bigl(\tilde{\theta}_j(g^1(x))-\tilde{\theta}_j(g^2(x))\Bigr)
\Bigl(\tilde{\theta}_n(g^1(x))-\tilde{\theta}_n(g^2(x))\Bigr) \, dx\,, \quad
A^q_{jm} =\int_\Omega \tilde{\chi}_j(x) \tilde{\chi}_m(x)\, dx\,, \\
&\mathfrak{r}^f_j = \int_\Omega \Bigl(\tilde{r}^{1,(k)}(x)-\tilde{r}^{2,(k)}(x)\Bigr) 
\Bigl(\tilde{\theta}_j(g^1(x))-\tilde{\theta}_j(g^2(x))\Bigr)\, dx\,,\\
&\mathfrak{r}^q_j = \int_\Omega \bigl(\tilde{r}^{1,(k)}(x) - \tilde{f}^{(k+1)}(g^1(x))\bigr)  \tilde{\chi}_j(x)\, dx \text{ or }\
\mathfrak{r}^q_j = \int_\Omega \bigl(\tilde{\tilde{r}}^{1,(k)}(x) - \tilde{f}^{(k+1)}(g^1(x))\bigr)  \tilde{\chi}_j(x)\, dx.
\end{aligned}
\]
For the least squares formulation
\[
\int_\Omega \sum_{i\in\{1,2\}}
\Bigl(\sum_{m=1}^M a_m\, \tilde{\chi}_m(x)+
\sum_{n=1}^N b_n\, \tilde{\theta}_n(g^i(x))-\tilde{r}^{i,(k)}(x)\Bigr)^2\, dx 
= \min_{\vec{a},\vec{b}}!
\]
of \eqref{fixedpoint_finaltime_overall} we solve the block matrix system
\[
\left(\begin{array}{cc}A^{qq}&A^{qf}\\(A^{qf})^T&A^{ff}\end{array}\right)
\left(\begin{array}{c}\vec{a}\\ \vec{b}\end{array}\right)
\left(\begin{array}{c}\mathfrak{r}^q\\ \mathfrak{r}^f\end{array}\right)
\]
with
\[
\begin{aligned}
&A^{qq}_{jm} =2\int_\Omega \tilde{\chi}_j(x) \tilde{\chi}_m(x)\, dx\,, \quad
A^{ff}_{jn} = \sum_{i\in\{1,2\}}\int_\Omega \tilde{\theta}_j(g^i(x))\, \tilde{\theta}_n(g^i(x)) \, dx\,, \\
&A^{qf}_{jn} = \sum_{i\in\{1,2\}}\int_\Omega \tilde{\chi}_j(x) \, \tilde{\theta}_n(g^i(x)) \, dx\,, 
\\
&\mathfrak{r}^q_j = \sum_{i\in\{1,2\}}\int_\Omega \tilde{r}^{i,(k)}(x) \, \tilde{\chi}_j(x)\, dx, \quad
\mathfrak{r}^f_j = \sum_{i\in\{1,2\}}\int_\Omega \tilde{r}^{i,(k)}(x) \, \tilde{\theta}_j(g^i(x))\, dx\, dx.
\end{aligned}
\]
for $i\in\{1,2\}$.

To take into account \eqref{qbar}, we consider in place of \eqref{basis_fixedpoint_finaltime}
\begin{equation}\label{basis_fixedpoint_finaltime_qbar}
\tilde{q}^{(k+1)}(x)=\log(\bar{q})+\sum_{m=1}^M a_M\, \tilde{\chi}_m(x), \qquad
\tilde{f}^{(k+1)}(u)=\sum_{n=1}^N b_n\, \tilde{\theta}_n(u)
\end{equation}
with basis functions $\tilde{\chi}_m$ that vanish at $\bar{x}$.
The system matrices remain as above, we only need to correct the right hand side vectors by substituting $\tilde{r}^{i,(k)}(x)$ by $\tilde{r}^{i,(k)}(x)-\log(\bar{q})$.

\subsubsection*{Fixed point scheme in case (c)}
To avoid division by zero at the left hand endpoint in \eqref{q_kp1} if $f(0)=0$, $g(0)=0$ (e.g., due to Dirichlet boundary conditions), we take $x_0>0$.

Making an ansatz
\[
f^{(k+1)}=\sum_{n=1}^N b_n\, \theta_n(u)
\]
we can reformulate \eqref{f_kp1} as a least squares problem
\[
\int_0^T \Bigl(\mathfrak{b}^{(k)}(t)-q^1\sum_{n=1}^N b_n\, \theta_n(h(t))\Bigr)^2\, dt = \min_{\vec{b}}!
\]
with $\mathfrak{b}^{(k)}(t)= -h'(t)+u^{(k)}_t(x_0,t) +q_1 f^{(k)}(h(t))$,
which leads to the linear system $A \vec{b}= \mathfrak{r}$ with 
\[
A_{jn} = q_1\int_0^T \theta_j(h(t)) \theta_n(h(t))\, dt\,, \quad
\mathfrak{r}_j = \int_0^T\mathfrak{b}^{(k)}(t) \theta_j(h(t))\, dt.
\]
Alternatively we can discretise the increment
\[
\df(u)=\sum_{n=1}^N b_n\, \theta_n(u)
\]
and appoximate the solution to \eqref{f_kp1} by setting $f^{(k+1)}=f^{(k)}+\df$, where $\vec{b}$ solves $A \vec{b}= \mathfrak{r}$ with $\mathfrak{r}_m = \int_0^T\mathfrak{b}^{(k)}(t) \theta_m(h(t))\, dt$, $\mathfrak{b}^{(k)}(t)= -h'(t)+u^{(k)}_t(x_0,t)$ and $A$ as above.

Instead of just dividing according to \eqref{q_kp1} (which does not enforce the boundary value $q(x_0)=q_1$ we also make an ansatz 
\[
q^{(k+1)}=\sum_{m=1}^M a_m\, \chi_m(x)
\]
where the basis functions $\chi_n$ vanish at $x_0$ and solve the least squares problem
\[
\int_\Omega \Bigl(\mathfrak{a}^{(k)}(x)-
\sum_{m=1}^M a_m\, \chi_m(x)\Bigr)^2\, dx = \min_{\vec{a}}!
\]
where 
$\mathfrak{a}^{(k)}(x)=\Bigl(-u^{(k)}_t(x,T)+\triangle g(x)+r(x,T)\Bigr)/f^{(k+1)}(g(x))$, which leads to the linear system $A \vec{a}= \mathfrak{r}$ with 
\[
A_{jn} = \int_\Omega \chi_j(x) \chi_n(x)\, dx\,, \quad
\mathfrak{r}_j = \int_\Omega \mathfrak{a}^{(k)}(x) \chi_j(x)\, dx.
\]